# Geometrically Intrinsic Nonlinear Recursive Filters II: Foundations

R. W. R. Darling [1]

[1] P. O. Box 535, Annapolis Junction, Maryland 20701-0535, USA. E-mail: rwrd@afterlife.ncsc.mil

**ABSTRACT:** This paper contains the technical foundations from stochastic differential geometry for the construction of geometrically intrinsic nonlinear recursive filters. A diffusion $X$ on a manifold $N$ is run for a time interval $T$, with a random initial condition. There is a single observation consisting of a nonlinear function of $X(T)$, corrupted by noise, and with values in another manifold $M$. The noise covariance of $X$ and the observation covariance themselves induce geometries on $M$ and $N$, respectively. Using these geometries we compute approximate but coördinate-free formulas for the "best estimate" of $X(T)$, given the observation, and its conditional variance. Calculations are based on use of Jacobi fields and of "intrinsic location parameters", a notion derived from the heat flow of harmonic mappings. When any nonlinearity is present, the resulting formulas are **not** the same as those for the continuous-discrete Extended Kalman Filter. A subsidiary result is a formula for computing approximately the "exponential barycenter" of a random variable $S$ on a manifold, i.e. a point $z$ such that the inverse image of $S$ under the exponential map at $z$ has mean zero in the tangent space at $z$.

**AMS (1991) SUBJECT CLASSIFICATION:** Primary: 60G35, 58G32. Secondary: 53B20

**KEY WORDS:** Jacobi field, stochastic differential equation, intrinsic location parameter, state estimation, nonlinear filter, extended Kalman filter

# 1    Overview

The purpose of this paper is to develop the necessary geometric and probabilistic machinery with which to construct the GI Filter algorithm presented in the companion paper, Darling [3]. The results can be summarized as follows.

## 1.1    State Process

Section 5.1 below explains how a diffusion process $X$, with $C^3$ coefficients, on a manifold $N$ induces, through its covariance structure, a semi-definite metric $\langle .|. \rangle$ on the cotangent bundle, a connection $\Gamma$ on the tangent bundle (actually $\Gamma$ may not be unique), and a (possibly degenerate) Laplace-Beltrami operator $\Delta$ such that the generator of $X$ can be written in intrinsic form as $L \equiv \xi + \frac{1}{2}\Delta$ for some vector field $\xi$.

Given $x_0 \in N$, let $\{x_t, 0 \le t \le \delta\}$ be the solution of the ODE associated with the vector field $\xi$, started at $x_0$. In other words, $x_t = \phi_t(x)$, for $0 \le t \le \delta$, where $\{\phi_t, 0 \le t \le \delta\}$ is the flow of the vector field $\xi$ on $N$.

Using the exponential map associated with $\Gamma$, we set up a random initial condition $X_0 \equiv \exp_{x_0}(U_0)$, where $U_0$ is a random variable in $T_{x_0}N$, with covariance $\Sigma_0 \in T_{x_0}N \otimes T_{x_0}N$. Then we run the diffusion during the time interval $[0, \delta]$, to obtain a **state process** $\{X_t, 0 \le t \le \delta\}$.

## 1.2    Observation

Now we take a nonlinear, noisy observation of the final value $X_\delta$. To be precise, we are given a $C^3$ function $\psi: N \to M$, where $M$ is a Riemannian manifold of dimension $q$. Let $\beta(y) \in T_yM \otimes T_yM$ be the inverse metric tensor at $y \in M$, which can be interpreted as the covariance of a random vector in $T_yM$. The **observation** $Y_1$ is of the form:

$$Y_1 \equiv \exp_{\psi(X_\delta)} V_1 \in M,$$

where $V_1$ is a mean-zero random vector in $T_{\psi(X_\delta)}M$, whose covariance is $\beta(y)$ when $\psi(X_\delta) = y$, but which is otherwise independent of $U_0$ and of the state process.

## 1.3    Intrinsic Approximation of the Conditional Distribution

Consider the random vector $U_\delta \oplus Z_\delta \in T_{x_0}N \oplus T_{\psi(x_\delta)}M$ given by

$$U_\delta \equiv \exp^{-1}_{x_\delta}(X_\delta), \ Z_\delta \equiv \exp^{-1}_{\psi(x_\delta)}(Y_1).$$

Our main result - Theorem 6.3 - presents intrinsic (i.e. coördinate-free) approximations to $\mathrm{E}[U_\delta|Z_\delta]$ and $\mathrm{Var}(U_\delta|Z_\delta)$, valid up to $O(\gamma^4)$, where $\Sigma_0$, the noise covariance of $X$, and the observation covariance are taken to be $O(\gamma^2)$, for some small $\gamma$, and $\|\mathrm{E}[U_0]\| = O(\gamma^4)$. Since Brownian scaling implies that $O(\gamma^2) = O(\delta)$, the remainder terms in our approximations are

$O(\delta^2)$. The form of the expression for $\text{Var}(U_\delta|Z_\delta)$ is quadratic in $Z_\delta$, unlike the non-intrinsic linear formula given in the Extended Kalman Filter.

### 1.4    Exponential Barycentres

In the light of the recursive nature of a filter, we want to use $Y_1$ to estimate the initial value for the diffusion $X_t' = X_{\delta+t}$, for $0 \le t \le \delta$, started at a point $X_0' \equiv \exp_{x_0'}(U_0')$ where

$$\|\text{E}[U_0']\| = O(\gamma^4), \|\text{E}[T(U_0' \otimes U_0' \otimes U_0')]\| = O(\gamma^4)$$

for tensor fields $T$ of type $(1, 3)$; in other words, $x_0'$ needs to be our best approximation to $X_\delta$, given $Y_1$. Our other major result is Theorem 3.2 (the Exponential Barycentre Formula), which implies that we should choose

$$x_0' \equiv \exp_{x_\delta}\{\hat{\mu} - \frac{1}{3}\sum_{i,j,k} R_{ijk}(x_\delta)\hat{\mu}^i\hat{\Sigma}^{jk}\}$$

where $\hat{\mu} \equiv \text{E}[U_\delta|Z_\delta]$, $\hat{\Sigma} \equiv \text{Var}(U_\delta|Z_\delta)$, $R_{ijk} \equiv R\left(\frac{\partial}{\partial x_i}, \frac{\partial}{\partial x_j}\right)\frac{\partial}{\partial x_k}$, and $R$ is the curvature tensor on $N$.

## 2    Curvature and the Exponential Map

The purpose of this section is to build some geometric machinery which will enable us to find asymptotic expansions for the first and second moments of a random variable of the form

$$\zeta \equiv \exp_y^{-1}(\exp_{\psi(X_\delta)}\eta) \in T_yM,$$

where $X$ is a diffusion process on a manifold $N$, $\psi: N \to M$ is a $C^2$ map, and $\eta$ is a random variable in $T_{\psi(X_\delta)}M$.

### 2.1    Geometric Notation

If $\{V(t), 0 \le t \le 1\}$ is a vector field along a $C^2$ path $\{y(t), 0 \le t \le 1\}$ on a manifold $M$ with a torsion-free connection, with local connector $\Gamma(.)$, then the covariant derivative of $V$ along $y$ is represented in coördinates by

$$\frac{\nabla V}{\partial t} \equiv V' + \Gamma(y)(V \otimes y'). \tag{1}$$

Recall that $V$ is said to be **parallel** along $y$ if (1) is identically zero, and $y$ is said to be a **geodesic** if $\nabla y'/\partial t$ is identically zero. For more information, see, for example, do Carmo [4] and Darling [1]. Let us recall the formula for the **curvature tensor** $R$ in terms of the local connector:

$$R(u,v)w = D\Gamma(v)(w \otimes u) - D\Gamma(u)(w \otimes v) + \Gamma(\Gamma(w \otimes u) \otimes v) - \Gamma(\Gamma(w \otimes v) \otimes u), \tag{2}$$

where $R, \Gamma, D\Gamma(v)$ are short for $R(y), \Gamma(y), D\Gamma(y)(v)$, respectively, and $y \in M$, $u, v, w \in T_yM$. This expression is alternating in $u$ and $v$.

## 2.2    Lemma

*Suppose W is a vector field along a geodesic $\{b(t), 0 \leq t \leq 1\}$ and $Z \equiv \frac{\nabla W}{\partial t}$. Suppose $b' = V$, and $\|V\|, \|W\|, \|Z\|$ are all $O(\gamma)$ at $t = 0$, for some small number $\gamma$. The non-intrinsic expansion in local coördinates for the tangent vector $W(1) \in T_{b(1)}M$, in terms of tangent vectors in $T_{b(0)}M$, is:*

$$W(1) = \left\{ W + \frac{\nabla W}{\partial t} + \frac{1}{2}\frac{\nabla^2 W}{\partial t^2} + \frac{1}{6}\frac{\nabla^3 W}{\partial t^3} + \frac{1}{2}R\left(V, W + \frac{\nabla W}{\partial t}\right)V \right\}\bigg|_{t=0} + \quad (3)$$

$$\left\{ -\Gamma\left(\left[W + Z + \frac{1}{2}\frac{\nabla Z}{\partial t}\right] \otimes V\right) + \Gamma(\Gamma(V \otimes [W+Z]) \otimes V) - \frac{1}{2}D\Gamma(W+Z)(V \otimes V) \right\}\bigg|_{t=0} + O(\gamma^4)$$

.

**Proof:** We see from (1) that

$$W' \equiv \frac{\partial W}{\partial t} = Z - \Gamma(b)(W \otimes b'),$$

$$\frac{\partial}{\partial t}\left(\frac{\nabla W}{\partial t}\right) = \frac{\nabla Z}{\partial t} - \Gamma(b)(Z \otimes b'),$$

$$\frac{\partial}{\partial t}\left(\frac{\nabla^2 W}{\partial t^2}\right) = \frac{\nabla^2 Z}{\partial t^2} - \Gamma(b)\left(\frac{\nabla Z}{\partial t} \otimes b'\right).$$

The second and third equations can be rewritten, by (1), as

$$\frac{\partial^2 W}{\partial t^2} = \frac{\nabla Z}{\partial t} - \Gamma(b)(Z \otimes b') - \frac{\partial}{\partial t}\{\Gamma(b)(W \otimes b')\}$$

$$\frac{\partial^3 W}{\partial t^3} = \frac{\nabla^2 Z}{\partial t^2} - \Gamma(b)\left(\frac{\nabla Z}{\partial t} \otimes b'\right) - \frac{\partial}{\partial t}\{\Gamma(b)(Z \otimes b')\} - \frac{\partial^2}{\partial t^2}\{\Gamma(b)(W \otimes b')\}$$

Using the geodesic equation $b'' + \Gamma(b)(b' \otimes b') = 0$,

$$\frac{\partial}{\partial t}\{\Gamma(b)(W \otimes b')\} = D\Gamma(b)(b')(W \otimes b') + \Gamma(b)(W' \otimes b') - \Gamma(b)(W \otimes \Gamma(b)(b' \otimes b'));$$

$$\frac{\partial^2}{\partial t^2}\{\Gamma(b)(W \otimes b')\} = 2D\Gamma(b')(W' \otimes b') + \Gamma(W'' \otimes b') - 2\Gamma(W' \otimes \Gamma(b' \otimes b')) + O(\gamma^4),$$

writing $\Gamma$ instead of $\Gamma(b)$ in the last equation. Take a Taylor expansion at $t = 0$, and write $\Gamma$ instead of $\Gamma(b(0))$, etc.:

$$W(1) = W_0 + \frac{\partial W}{\partial t}(0) + \frac{1}{2}\frac{\partial^2 W}{\partial t^2}(0) + \frac{1}{6}\frac{\partial^3 W}{\partial t^3}(0) + O(\gamma^4)$$

$$= \{W + Z - \Gamma(W \otimes V)\}|_{t=0} + O(\gamma^4) +$$

$$\frac{1}{2}\{\frac{\nabla Z}{\partial t} - \Gamma(Z \otimes V) - D\Gamma(V)(W \otimes V) - \Gamma([Z - \Gamma(W \otimes V)] \otimes V) + \Gamma(W \otimes \Gamma(V \otimes V))\}\Big|_{t=0} +$$

$$\frac{1}{6}\left\{\frac{\nabla^2 Z}{\partial t^2} - \Gamma\left(\frac{\nabla Z}{\partial t} \otimes V\right) - D\Gamma(V)(Z \otimes V) - \Gamma\left(\left[\frac{\nabla Z}{\partial t} - \Gamma(Z \otimes V)\right] \otimes V\right) + \Gamma(Z \otimes \Gamma(V \otimes V))\right\}\Big|_{t=0}$$

$$+ \frac{1}{6}\{2\Gamma(Z \otimes \Gamma(V \otimes V)) - 2D\Gamma(V)(Z \otimes V) - \Gamma\left(\left[\frac{\nabla Z}{\partial t} - \Gamma(Z \otimes V) - \Gamma(Z \otimes V)\right] \otimes V\right)\}\Big|_{t=0}.$$

The terms involving one vector field are as in (3). Terms involving two vector fields are

$$-\Gamma(W \otimes V) - \Gamma(Z \otimes V) - \frac{1}{2}\Gamma\left(\frac{\nabla Z}{\partial t} \otimes V\right) = -\Gamma\left(\left[W + Z + \frac{1}{2}\frac{\nabla Z}{\partial t}\right] \otimes V\right).$$

Terms involving three vector fields are

$$\frac{1}{2}D\Gamma(V)(W \otimes V) + \frac{1}{2}\Gamma(\Gamma(W \otimes V) \otimes V) + \frac{1}{2}\Gamma(W \otimes \Gamma(V \otimes V))$$

$$+ \frac{1}{6}\{3\Gamma(\Gamma(Z \otimes V) \otimes V) + 3\Gamma(Z \otimes \Gamma(V \otimes V)) - 3D\Gamma(V)(Z \otimes V)\}$$

$$= -\frac{1}{2}D\Gamma(V)(V \otimes [W + Z]) + \frac{1}{2}\Gamma(\Gamma(V \otimes V) \otimes [W + Z]) + \frac{1}{2}\Gamma(\Gamma(V \otimes [W + Z]) \otimes V)$$

$$= \frac{1}{2}R(V, W + Z)V + \Gamma(\Gamma(V \otimes [W + Z]) \otimes V) - \frac{1}{2}D\Gamma(W + Z)(V \otimes V).$$

All terms in (3) have now been obtained. ◊

### 2.3   Construction

If $V$ is a vector field along a $C^2$ path $\{y(\varepsilon), 0 \le \varepsilon \le 1\}$ on a Riemannian manifold $M$, define a function $b: [0, 1] \times [0, 1] \to M$ by

$$b_\varepsilon(t) \equiv b(t, \varepsilon) \equiv \exp_{y(\varepsilon)} tV(\varepsilon). \tag{4}$$

Thus for each $\varepsilon$, $b_\varepsilon$ is a geodesic with

$$b_\varepsilon(0) = y(\varepsilon), \ b_\varepsilon'(0) = V(\varepsilon). \tag{5}$$

For each $\varepsilon \in [0, 1]$,

$$J_\varepsilon(t) \equiv \frac{\partial b}{\partial \varepsilon}(t, \varepsilon).$$

is a "Jacobi field" along the geodesic $b_\varepsilon$ with

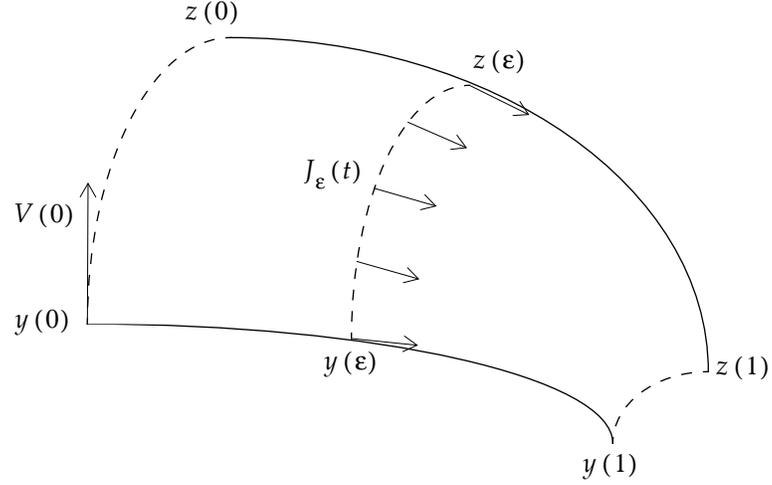

$$J_\varepsilon(0) = y'(\varepsilon), \quad \frac{\nabla J_\varepsilon}{\partial t}(0) = \frac{\nabla V}{\partial \varepsilon}. \tag{6}$$

The last identity follows from a well-known formula (see do Carmo [4], p. 68):

$$\frac{\nabla}{\partial \varepsilon}\left(\frac{\partial b_\varepsilon}{\partial t}\right) = \frac{\partial^2 b_\varepsilon}{\partial \varepsilon \partial t} + \Gamma(b_\varepsilon)(b_\varepsilon' \otimes J_\varepsilon) = \frac{\nabla J_\varepsilon}{\partial t}, \tag{7}$$

where $b_\varepsilon'$ is the vector field $\partial b_\varepsilon / \partial t$ along $b_\varepsilon$. From do Carmo [4], p. 98, we have that for every vector field $W_\varepsilon$ along $b_\varepsilon$,

$$\frac{\nabla}{\partial t}\frac{\nabla W_\varepsilon}{\partial \varepsilon} = \frac{\nabla}{\partial \varepsilon}\frac{\nabla W_\varepsilon}{\partial t} + R(J_\varepsilon, b_\varepsilon') W_\varepsilon. \tag{8}$$

Taking $W_\varepsilon \equiv b_\varepsilon'$ and using the fact that $b_\varepsilon$ is a geodesic (i.e. $\dfrac{\nabla b_\varepsilon'}{\partial t} = 0$), we obtain from (7) and (8) the Jacobi equation:

$$\frac{\nabla^2 J_\varepsilon}{\partial t^2} + R(b_\varepsilon', J_\varepsilon) b_\varepsilon' = 0. \tag{9}$$

Moreover $J_\varepsilon(1) = z'(\varepsilon)$, where $z(\varepsilon) \equiv \exp_{y(\varepsilon)} V(\varepsilon) \in M$.

In the following Proposition, only the case where $V(0) = 0$ will be needed in this paper; however the general case may be needed for other filtering investigations.

### 2.4　　Proposition

*Suppose $\{V(\varepsilon), 0 \leq \varepsilon \leq 1\}$ is a vector field along a $C^2$ path $\{y(\varepsilon), 0 \leq \varepsilon \leq 1\}$ on a Riemannian manifold M, such that, $\|y'\|, \|V\|, \left\|\dfrac{\nabla V}{\partial \varepsilon}\right\|$, are all $O(\gamma)$, for some small number $\gamma$. Define*

$$\zeta(\varepsilon) \equiv \exp_{y(0)}^{-1} (\exp_{y(\varepsilon)} V(\varepsilon)) \in T_{y(0)} M. \tag{10}$$

Then we have the following intrinsic formulae, whose local expressions follow from (1) and (2):

$$\zeta'(0) = \left\{ y' + \frac{\nabla V}{\partial \varepsilon} - \frac{1}{3} R(V, y') V \right\} \bigg|_{\varepsilon=0} + O(\gamma^4); \tag{11}$$

$$\zeta''(0) = \left\{ \frac{\nabla y'}{\partial \varepsilon} + \frac{\nabla^2 V}{\partial \varepsilon^2} - \frac{1}{3} R\left(V, \frac{\nabla y'}{\partial \varepsilon}\right) V - \frac{2}{3} R\left(\frac{\nabla V}{\partial \varepsilon}, y'\right) V + \frac{1}{3} R(y', V)\left(y' + 2\frac{\nabla V}{\partial \varepsilon}\right) \right\} \bigg|_{\varepsilon=0} + O(\gamma^4) \tag{12}$$

$$\zeta^{(3)}(0) = \left\{ \frac{\nabla^2 y'}{\partial \varepsilon^2} + \frac{\nabla^3 V}{\partial \varepsilon^3} + R\left(y', \frac{\nabla V}{\partial \varepsilon}\right)\left(y' + 2\frac{\nabla V}{\partial \varepsilon}\right) + R(y', V)\left(\frac{\nabla y'}{\partial \varepsilon} + 2\frac{\nabla^2 V}{\partial \varepsilon^2}\right) \right\} \bigg|_{\varepsilon=0} \tag{13}$$

$$+ \left\{ R\left(V, \frac{\nabla^2 V}{\partial \varepsilon^2}\right) y' + R\left(V, \frac{\nabla V}{\partial \varepsilon}\right)\frac{\nabla y'}{\partial \varepsilon} + 2R\left(\frac{\nabla y'}{\partial \varepsilon}, V\right)\frac{\nabla V}{\partial \varepsilon} - \frac{1}{3} R\left(V, \frac{\nabla^2 y'}{\partial \varepsilon^2}\right) V \right\} \bigg|_{\varepsilon=0} + O(\gamma^4).$$

**Proof: Step I.** Fix $\varepsilon$, and apply Lemma 2.2 with $W = J_\varepsilon$, using the initial conditions (6) and (9). Note that (9) and the geodesic equation imply that

$$\frac{\nabla^3 J_\varepsilon}{\partial t^3} = -R\left(b_\varepsilon', \frac{\nabla J_\varepsilon}{\partial t}\right) b_\varepsilon' + O(\gamma^4).$$

Consequently

$$\left\{ W + \frac{\nabla W}{\partial t} + \frac{1}{2}\frac{\nabla^2 W}{\partial t^2} + \frac{1}{6}\frac{\nabla^3 W}{\partial t^3} + \frac{1}{2} R\left(V, W + \frac{\nabla W}{\partial t}\right) V \right\} \bigg|_{t=0} =$$

$$y' + \frac{\nabla V}{\partial \varepsilon} - \frac{1}{2} R(V, y') V - \frac{1}{6} R\left(V, \frac{\nabla V}{\partial \varepsilon}\right) V + \frac{1}{2} R\left(V, y' + \frac{\nabla V}{\partial \varepsilon}\right) V + O(\gamma^4)$$

$$= y' + \frac{\nabla V}{\partial \varepsilon} + \frac{1}{3} R\left(V, \frac{\nabla V}{\partial \varepsilon}\right) V + O(\gamma^4). \tag{14}$$

The other terms in (3) contribute

$$-\Gamma\left(\left[y' + \frac{\nabla V}{\partial \varepsilon}\right] \otimes V\right) + \Gamma\left(\Gamma\left(V \otimes \left[y' + \frac{\nabla V}{\partial \varepsilon}\right]\right) \otimes V\right) - \frac{1}{2} D\Gamma\left(y' + \frac{\nabla V}{\partial \varepsilon}\right)(V \otimes V) + O(\gamma^4) \tag{15}$$

Thus $J_\varepsilon(1) = z'(\varepsilon)$ is approximated by the sum of (14) and (15).

**Step II.** Now apply Lemma 2.2 again, this time to the Jacobi field $K_\varepsilon$ along

$$c_\varepsilon(t) \equiv \exp_{y(0)} t\zeta(\varepsilon),$$

with $\zeta(\varepsilon)$ in place of $V(\varepsilon)$ and $y(0)$ in place of $y(\varepsilon)$. We obtain expressions analogous to (14) and (15), namely:

$$K_\varepsilon(1) = \zeta' + \frac{1}{3}R(\zeta, \zeta')\zeta - \Gamma(\zeta' \otimes V) + \Gamma(\Gamma(V \otimes \zeta') \otimes V) - \frac{1}{2}D\Gamma(\zeta')(V \otimes V) + O(\gamma^4). \quad (16)$$

Since, by definition,

$$z(\varepsilon) \equiv \exp_{y(\varepsilon)} V(\varepsilon) = \exp_{y(0)} \zeta(\varepsilon),$$

it follows that

$$J_\varepsilon(1) = z'(\varepsilon) = K_\varepsilon(1). \quad (17)$$

At $\varepsilon = 0$, $\zeta(0) = V(0)$. Let us adopt the working hypothesis that

$$\zeta'(0) = \{y' + \frac{\nabla V}{\partial \varepsilon} + S\}\bigg|_{\varepsilon = 0} + O(\gamma^4),$$

where $S$ is the sum of terms of third order. Combining (14), (15), and (16), we find

$$y' + \frac{\nabla V}{\partial \varepsilon} + \frac{1}{3}R\left(V, \frac{\nabla V}{\partial \varepsilon}\right)V = y' + \frac{\nabla V}{\partial \varepsilon} + S + \frac{1}{3}R\left(V, y' + \frac{\nabla V}{\partial \varepsilon}\right)V + O(\gamma^4),$$

giving $S = -\frac{1}{3}R(V, y')V$. This verifies (11).

**Step III.** Apply Lemma 2.2 a third time, taking $W = \frac{\nabla J_\varepsilon}{\partial \varepsilon}$. Using (7) and (8),

$$\frac{\nabla W}{\partial t} = \frac{\nabla}{\partial t}\left(\frac{\nabla J_\varepsilon}{\partial \varepsilon}\right) = \frac{\nabla}{\partial \varepsilon}\frac{\nabla J_\varepsilon}{\partial t} + R(J_\varepsilon, b_\varepsilon')J_\varepsilon; \quad (18)$$

$$\frac{\nabla^2 W}{\partial t^2} = \frac{\nabla}{\partial t}\{\frac{\nabla}{\partial t}\left(\frac{\nabla J_\varepsilon}{\partial \varepsilon}\right)\} = \frac{\nabla}{\partial t}\{\frac{\nabla}{\partial \varepsilon}\frac{\nabla J_\varepsilon}{\partial t} + R(J_\varepsilon, b_\varepsilon')J_\varepsilon\}$$

$$= \frac{\nabla}{\partial \varepsilon}\left\{\frac{\nabla^2 J_\varepsilon}{\partial t^2}\right\} + R(J_\varepsilon, b_\varepsilon')\frac{\nabla J_\varepsilon}{\partial t} + R\left(\frac{\nabla J_\varepsilon}{\partial t}, b_\varepsilon'\right)J_\varepsilon + R(J_\varepsilon, b_\varepsilon')\frac{\nabla J_\varepsilon}{\partial t} + O(\gamma^4)$$

$$= \frac{\nabla}{\partial \varepsilon}\{-R(b_\varepsilon', J_\varepsilon)b_\varepsilon'\} + 2R(J_\varepsilon, b_\varepsilon')\frac{\nabla J_\varepsilon}{\partial t} + R\left(\frac{\nabla J_\varepsilon}{\partial t}, b_\varepsilon'\right)J_\varepsilon + O(\gamma^4)$$

$$= -R\left(\frac{\nabla J_\varepsilon}{\partial t}, J_\varepsilon\right)b_\varepsilon' - R\left(b_\varepsilon', \frac{\nabla J_\varepsilon}{\partial \varepsilon}\right)b_\varepsilon' - R(b_\varepsilon', J_\varepsilon)\frac{\nabla J_\varepsilon}{\partial t} + 2R(J_\varepsilon, b_\varepsilon')\frac{\nabla J_\varepsilon}{\partial t} + R\left(\frac{\nabla J_\varepsilon}{\partial t}, b_\varepsilon'\right)J_\varepsilon + O(\gamma^4)$$

$$= -R\left(\frac{\nabla J_\varepsilon}{\partial t}, J_\varepsilon\right)b_\varepsilon' - R\left(b_\varepsilon', \frac{\nabla J_\varepsilon}{\partial \varepsilon}\right)b_\varepsilon' + 3R(J_\varepsilon, b_\varepsilon')\frac{\nabla J_\varepsilon}{\partial t} - R\left(b_\varepsilon', \frac{\nabla J_\varepsilon}{\partial \varepsilon}\right)J_\varepsilon + O(\gamma^4).$$

We used $\frac{\nabla b_\varepsilon'}{\partial t} = 0$ here. Using the first Bianchi identity,

$$R\left(\frac{\nabla J_\varepsilon}{\partial t}, J_\varepsilon\right)b_\varepsilon' + R(J_\varepsilon, b_\varepsilon')\frac{\nabla J_\varepsilon}{\partial t} + R\left(b_\varepsilon', \frac{\nabla J_\varepsilon}{\partial \varepsilon}\right)J_\varepsilon = 0,$$

$$\frac{\nabla^2 W}{\partial t^2} = -R\left(b_\varepsilon', \frac{\nabla J_\varepsilon}{\partial \varepsilon}\right) b_\varepsilon' + 4R(J_\varepsilon, b_\varepsilon') \frac{\nabla J_\varepsilon}{\partial t} + O(\gamma^4). \tag{19}$$

Finally, using (18) and (19),

$$\frac{\nabla^3 W}{\partial t^3} = \frac{\nabla}{\partial t}\{-R\left(b_\varepsilon', \frac{\nabla J_\varepsilon}{\partial \varepsilon}\right) b_\varepsilon' + 4R(J_\varepsilon, b_\varepsilon') \frac{\nabla J_\varepsilon}{\partial t}\} + O(\gamma^4)$$

$$= -R\left(b_\varepsilon', \frac{\nabla}{\partial t}\left(\frac{\nabla J_\varepsilon}{\partial \varepsilon}\right)\right) b_\varepsilon' + 4R\left(\frac{\nabla J_\varepsilon}{\partial t}, b_\varepsilon'\right) \frac{\nabla J_\varepsilon}{\partial t} + O(\gamma^4)$$

$$= -R\left(b_\varepsilon', \frac{\nabla}{\partial \varepsilon}\frac{\nabla J_\varepsilon}{\partial t}\right) b_\varepsilon' + 4R\left(\frac{\nabla J_\varepsilon}{\partial t}, b_\varepsilon'\right) \frac{\nabla J_\varepsilon}{\partial t} + O(\gamma^4).$$

Thus our initial conditions will be

$$W(0) = \frac{\nabla y'}{\partial \varepsilon};$$

$$\frac{\nabla W}{\partial t}(0) = \frac{\nabla^2 V}{\partial \varepsilon^2} + R(y', V) y';$$

$$\frac{\nabla^2 W}{\partial t^2}(0) = -R\left(V, \frac{\nabla y'}{\partial \varepsilon}\right) V + 4R(y', V) \frac{\nabla V}{\partial \varepsilon} + O(\gamma^4);$$

$$\frac{\nabla^3 W}{\partial t^3}(0) = -R\left(V, \frac{\nabla^2 V}{\partial \varepsilon^2}\right) V + 4R\left(\frac{\nabla V}{\partial \varepsilon}, V\right) \frac{\nabla V}{\partial \varepsilon} + O(\gamma^4).$$

By Lemma 2.2,

$$\frac{\nabla J_\varepsilon}{\partial \varepsilon}(1) = \frac{\nabla y'}{\partial \varepsilon} + \frac{\nabla^2 V}{\partial \varepsilon^2} + R(y', V) y' + \frac{1}{2}\{-R\left(V, \frac{\nabla y'}{\partial \varepsilon}\right) V + 4R(y', V) \frac{\nabla V}{\partial \varepsilon}\} + \tag{20}$$

$$\frac{1}{6}\left\{-R\left(V, \frac{\nabla^2 V}{\partial \varepsilon^2}\right) V + 4R\left(\frac{\nabla V}{\partial \varepsilon}, V\right) \frac{\nabla V}{\partial \varepsilon}\right\} + \frac{1}{2}R\left(V, \frac{\nabla y'}{\partial \varepsilon} + \frac{\nabla^2 V}{\partial \varepsilon^2}\right) V + O(\gamma^4) +$$

$$-\Gamma\left(\left[\frac{\nabla y'}{\partial \varepsilon} + \frac{\nabla^2 V}{\partial \varepsilon^2}\right] \otimes V\right) + \Gamma\left(\Gamma\left(V \otimes \left[\frac{\nabla y'}{\partial \varepsilon} + \frac{\nabla^2 V}{\partial \varepsilon^2}\right]\right) \otimes V\right) - \frac{1}{2}D\Gamma\left(\frac{\nabla y'}{\partial \varepsilon} + \frac{\nabla^2 V}{\partial \varepsilon^2}\right)(V \otimes V).$$

The curvature terms simplify to:

$$R(y', V) y' + \frac{1}{3}R\left(V, \frac{\nabla^2 V}{\partial \varepsilon^2}\right) V + \frac{2}{3}R\left(\frac{\nabla V}{\partial \varepsilon}, V\right) \frac{\nabla V}{\partial \varepsilon} + 2R(y', V) \frac{\nabla V}{\partial \varepsilon}.$$

**Step IV.** With reference to Step II, apply Lemma 2.2 a fourth time, taking $W = \dfrac{\nabla K_\varepsilon}{\partial \varepsilon}$. The expression analogous to (20) is:

$$\frac{\nabla K_\varepsilon}{\partial \varepsilon}(1) = \zeta'' + \frac{1}{3}R(\zeta, \zeta'')\zeta + \frac{2}{3}R(\zeta', \zeta)\zeta'$$

$$-\Gamma(\zeta'' \otimes \zeta) + \Gamma(\Gamma(\zeta \otimes \zeta'') \otimes V) - \frac{1}{2}D\Gamma(\zeta'')(\zeta \otimes \zeta) + O(\gamma^4).$$

From (17), it follows that $\dfrac{\nabla J_\varepsilon}{\partial \varepsilon}(1) = \dfrac{\nabla K_\varepsilon}{\partial \varepsilon}(1)$. Take as a working hypothesis that

$$\zeta''(0) = \left\{ \frac{\nabla y'}{\partial \varepsilon} + \frac{\nabla^2 V}{\partial \varepsilon^2} + T \right\}\bigg|_{\varepsilon=0} + O(\gamma^4) \tag{21}$$

where $T$ is a sum of third order (i.e. curvature) terms. Combining (11) and (20) - (21) gives

$$T + \frac{1}{3}R\left(V, \frac{\nabla y'}{\partial \varepsilon} + \frac{\nabla^2 V}{\partial \varepsilon^2}\right)V + \frac{2}{3}R\left(y' + \frac{\nabla V}{\partial \varepsilon}, V\right)\left(y' + \frac{\nabla V}{\partial \varepsilon}\right) =$$

$$R(y', V)y' + \frac{1}{3}R\left(V, \frac{\nabla^2 V}{\partial \varepsilon^2}\right)V + \frac{2}{3}R\left(\frac{\nabla V}{\partial \varepsilon}, V\right)\frac{\nabla V}{\partial \varepsilon} + 2R(y', V)\frac{\nabla V}{\partial \varepsilon};$$

$$T + \frac{1}{3}R\left(V, \frac{\nabla y'}{\partial \varepsilon}\right)V + \frac{2}{3}R\left(\frac{\nabla V}{\partial \varepsilon}, V\right)y' = \frac{1}{3}R(y', V)y' + \frac{4}{3}R(y', V)\frac{\nabla V}{\partial \varepsilon};$$

$$T = \frac{1}{3}R(y', V)y' + \frac{4}{3}R(y', V)\frac{\nabla V}{\partial \varepsilon} + \frac{2}{3}R\left(V, \frac{\nabla V}{\partial \varepsilon}\right)y' - \frac{1}{3}R\left(V, \frac{\nabla y'}{\partial \varepsilon}\right)V.$$

Using the first Bianchi identity,

$$R\left(V, \frac{\nabla V}{\partial \varepsilon}\right)y' + R\left(\frac{\nabla V}{\partial \varepsilon}, y'\right)V + R(y', V)\frac{\nabla V}{\partial \varepsilon} = 0; \tag{22}$$

$$T = \frac{1}{3}R(y', V)y' + \frac{2}{3}R(y', V)\frac{\nabla V}{\partial \varepsilon} - \frac{2}{3}R\left(\frac{\nabla V}{\partial \varepsilon}, y'\right)V - \frac{1}{3}R\left(V, \frac{\nabla y'}{\partial \varepsilon}\right)V.$$

Thus (12) is established. The proof of (13) is similar, taking $W = \dfrac{\nabla^2 J_\varepsilon}{\partial \varepsilon^2}$.  ◊

We conclude this section with two minor formulas we will need to implement the filter.

### 2.5     Lemma (Expansion of the Exponential Map in Local Coordinates)

*On a manifold M with a torsion-free connection $\Gamma$, the expansion of the exponential map in local coordinates is given in terms of the local connector $\Gamma(.)$ as follows: for $y \in M \cong R^q$, $v \in T_y M \cong R^q$, $t \in R$,*

$$\exp_y tv =$$

$$y + tv - \frac{t^2}{2}\Gamma(y)(v \otimes v) + \frac{t^3}{6}[2\Gamma(y)(\Gamma(y)(v \otimes v) \otimes v) - D\Gamma(y)(v)(v \otimes v)] + O(t^4). \tag{23}$$

**Proof:** The geodesic $b(t) \equiv \exp_y tv$ has the following properties. By definition $b(0) = y$ and $b'(0) = v$, and by the geodesic equation, $b'' + \Gamma(b)(b' \otimes b') = 0$, which implies that $b''(0)$ is $-\Gamma(y)(v \otimes v)$. A further differentiation of the geodesic equation yields

$$b^{(3)} + D\Gamma(b)(b')(b' \otimes b') + \Gamma(b)(b'' \otimes b') + \Gamma(b)(b' \otimes b'') = 0,$$

and the last two terms on the left are the same by torsion-freedom. Now set $t = 0$ and substitute for $b'(0)$ and $b''(0)$ to give the expression for $b^{(3)}(0)$ which appears as the coefficient of $t^3/6$ in the Taylor expansion (23). ◊

## 2.6   Corollary (Expansion of the Inverse Exponential Map)

*The corresponding expansion for $\exp_y^{-1}(z)$, taking $w \equiv z - y$, is:*

$$w + \frac{1}{2}\Gamma(y)(w \otimes w) + \frac{1}{6}[D\Gamma(y)(w)(w \otimes w) + \Gamma(y)(\Gamma(y)(w \otimes w) \otimes w)] + O(\|w\|^4). \qquad (24)$$

**Proof:** Let $b(t) \equiv \exp_y tv$ and $h(z) \equiv \exp_y^{-1}(z)$, so that $h \bullet b(t) = tv$. Applying the chain rule repeatedly,

$$Dh(b(t)) \bullet b'(t) = v,$$

$$D^2 h(b(t))(b' \otimes b') + Dh(b(t)) \bullet b''(t) = 0,$$

$$D^3 h(b(t))(b' \otimes b' \otimes b') + 3D^2 h(b(t))(b'' \otimes b') + Dh(b(t)) \bullet b^{(3)}(t) = 0.$$

Now set $t = 0$, and substitute the expressions computed in the proof of Lemma 2.5:

$$Dh(y)v = v;$$

$$D^2 h(y)(v \otimes v) - \Gamma(y)(v \otimes v) = 0;$$

$$D^3 h(y)(v \otimes v \otimes v)$$
$$= 3\Gamma(y)(\Gamma(y)(v \otimes v) \otimes v) - [2\Gamma(y)(\Gamma(y)(v \otimes v) \otimes v) - D\Gamma(y)(v)(v \otimes v)].$$

Since $h(y) = 0$, we obtain (24) as the third order Taylor expansion of $h$ at $z$. ◊

## 3   Application to Computation of Exponential Barycentres

Let us recall from Emery and Mokobodzki [5] that an **exponential barycentre** for a random variable $Z$ on a manifold $N$ with a torsion-free connection $\Gamma$ is a point $z \in N$ such that the random variable $\exp_z^{-1}(Z) \in T_z N$ has mean zero.

### 3.1   Problem

The solution of the following problem is needed in order to update the GI Filter in Darling [3]. We are given a point $x \in N$, and moments

$$\mu \equiv E\,[\exp_x^{-1}(Z)] \in T_xN, \; \Sigma \equiv \mathrm{Var}\,(\exp_x^{-1}(Z)) \in T_xN \otimes T_xN. \tag{25}$$

We would like to compute from these moments an exponential barycentre for $Z$, or at least a "good" approximation.

### 3.1.a    Notation

In the following theorem, the norm $\|\,.\,\|$ is with respect to some reference metric for $N$, which need not be related to the connection. Given the curvature tensor (2), the vector field $R\!\left(\dfrac{\partial}{\partial x_i}, \dfrac{\partial}{\partial x_j}\right)\!\dfrac{\partial}{\partial x_k}$ is denoted $R_{ijk}$.

### 3.2    Theorem (Exponential Barycentre Formula)

*Suppose that the data (25) satisfy: $\|\mu\| = O(\gamma)$, $\|\Sigma\| = O(\gamma^2)$, for some small number $\gamma$. Define*

$$z \equiv \exp_x \{\mu - \frac{1}{3} \sum_{i,j,k} R_{ijk} \mu^i \Sigma^{jk}\}. \tag{26}$$

*Then*

$$E\,[\exp_z^{-1}(Z)] = O(\gamma^4) \tag{27}$$

*In other words, $z$ is an approximate exponential barycentre for $Z$. Also if $T$ is a tensor field of type $(1,3)$ (such as a curvature tensor), and if $\|E\,[T(\eta - \mu, \eta - \mu, \eta - \mu)]\| = O(\gamma^4)$, where $\eta \equiv \exp_x^{-1}(Z) \in T_xN$, then*

$$\left\| E\,[T(\exp_z^{-1}(Z), \exp_z^{-1}(Z), \exp_z^{-1}(Z))] \right\| = O(\gamma^4).$$

**Remark**: Lemma 2.5 gives a formula to compute the exponential barycentre $z$ in local coördinates.

**Proof:** Let $\{y(\varepsilon), 0 \leq \varepsilon \leq 1\}$ be a geodesic on $N$ with initial conditions:

$$y(0) \equiv x, \; y'(0) \equiv \mu - \rho,$$

where $\rho \equiv \dfrac{1}{3}\sum_{i,j,k} R_{ijk} \mu^i \Sigma^{jk} \in T_xN$, so that $\rho = O(\gamma^3)$.

Let $\eta$ denote the random variable $\exp_x^{-1}(Z) \in T_xN$, whose first and second moments are given by (25), and define a random vector field $V$ along $y$ by the formula

$$V(\varepsilon) \equiv \exp_{y(\varepsilon)}^{-1} \{\exp_x(\varepsilon \eta)\}.$$

Now $y(1) = z$ by (26), and $V(1) = \exp_z^{-1}(Z)$, so to prove (27) it suffices to prove that

$$E\,[V(1)] = O(\gamma^4). \tag{28}$$

We may rewrite the definition in the form (cf. (10))

$$\zeta(\varepsilon) \equiv \varepsilon\eta = \exp_{y(0)}^{-1}\{\exp_{y(\varepsilon)} V(\varepsilon)\}.$$

Now apply Proposition 2.4, noting that $V(0) = 0$, and since $y$ is a geodesic, $\frac{\nabla y'}{\partial \varepsilon} = \frac{\nabla^2 y'}{\partial \varepsilon^2} = 0$. We obtain:

$$\frac{\nabla V}{\partial \varepsilon}(0) = \eta - y'(0) + O(\gamma^4); \quad \frac{\nabla^2 V}{\partial \varepsilon^2}(0) = O(\gamma^4);$$

$$\frac{\nabla^3 V}{\partial \varepsilon^3}(0) = \left\{-R\left(y', \frac{\nabla V}{\partial \varepsilon}\right)y' - 2R\left(y', \frac{\nabla V}{\partial \varepsilon}\right)\frac{\nabla V}{\partial \varepsilon}\right\}\bigg|_{\varepsilon=0} + O(\gamma^4).$$

By Lemma 2.2, with $V$ in place of $W$,

$$V(1) = \left\{\frac{\nabla V}{\partial \varepsilon} + \frac{1}{6}\frac{\nabla^3 V}{\partial \varepsilon^3} + \frac{1}{2}R\left(y', \frac{\nabla V}{\partial \varepsilon}\right)y'\right\}\bigg|_{\varepsilon=0} +$$

$$\left\{-\Gamma\left(\frac{\nabla V}{\partial \varepsilon} \otimes y'\right) + \Gamma\left(\Gamma\left(y' \otimes \frac{\nabla V}{\partial \varepsilon}\right) \otimes y'\right) - \frac{1}{2}D\Gamma\left(\frac{\nabla V}{\partial \varepsilon}\right)(y' \otimes y')\right\}\bigg|_{\varepsilon=0} + O(\gamma^4).$$

Taking expectations, and substituting,

$$E[V(1)] = E\left[\eta - \mu + \rho + \frac{1}{6}\{-R(\mu, \eta)\mu - 2R(\mu, \eta)(\eta - \mu)\} + \frac{1}{2}R(\mu, \eta)\mu\right] +$$

$$E\left[-\Gamma((\eta - \mu) \otimes \mu) + \Gamma(\Gamma(\mu \otimes (\eta - \mu)) \otimes \mu) - \frac{1}{2}D\Gamma(\eta - \mu)(\mu \otimes \mu)\right] + O(\gamma^4).$$

The expected values in the second line vanish by linearity and the fact that $E[\eta] = \mu$. Also $E[R(\mu, \eta)\mu] = R(\mu, \mu)\mu = 0$. Now the terms on the right cancel up to $O(\gamma^4)$, since

$$E[R(\mu, \eta)\eta] = \sum_{i,j,k} R_{ijk}{}^i \mu^i E[\eta^j \eta^k] = 3\rho.$$

Thus (28) is established. As for the assertion about the tensor $T$, let $W$ be the random vector field along $y$ given by $W \equiv T(V, V, V)$. Then $W(0) = 0$,
$W(1) = T(\exp_z^{-1}(Z), \exp_z^{-1}(Z), \exp_z^{-1}(Z))$, and

$$\frac{\nabla W}{\partial \varepsilon}(0) = \left\{T\left(\frac{\nabla V}{\partial \varepsilon}, V, V\right) + T\left(V, \frac{\nabla V}{\partial \varepsilon}, V\right) + T\left(V, V, \frac{\nabla V}{\partial \varepsilon}\right)\right\}\bigg|_{\varepsilon=0} + O(\gamma^4) = O(\gamma^4)$$

$$\frac{\nabla^2 W}{\partial \varepsilon^2}(0) = \left\{2T\left(\frac{\nabla V}{\partial \varepsilon}, \frac{\nabla V}{\partial \varepsilon}, V\right) + 2T\left(V, \frac{\nabla V}{\partial \varepsilon}, \frac{\nabla V}{\partial \varepsilon}\right) + 2T\left(V, \frac{\nabla V}{\partial \varepsilon}, \frac{\nabla V}{\partial \varepsilon}\right)\right\}\bigg|_{\varepsilon=0} + O(\gamma^4) = O(\gamma^4)$$

$$\frac{\nabla^3 W}{\partial \varepsilon^3}(0) = 6\left\{T\left(\frac{\nabla V}{\partial \varepsilon}, \frac{\nabla V}{\partial \varepsilon}, \frac{\nabla V}{\partial \varepsilon}\right)\right\}\bigg|_{\varepsilon=0} + O(\gamma^4) = 6\{T(\eta - \mu, \eta - \mu, \eta - \mu)\}|_{\varepsilon=0} + O(\gamma^4)$$

where we have used the fact that $V(0) = 0$ and $\frac{\nabla^2 V}{\partial \varepsilon^2}(0) = O(\gamma^4)$. Now apply Lemma 2.2 and the assumption in the Theorem to deduce

$$\|E[W(1)]\| = O(\gamma^4)$$

as desired.                                                                                           ◊

# 4    A Non-Gaussian Conditional Expectation Formula

## 4.1    Quadratic Forms in Multivariate Normal Random Vectors

Suppose that $U$, $V$, and $Z$ are random vectors in $R^p$, $R^q$, and $R^r$, respectively, whose joint distribution is multivariate Normal with

$$\begin{bmatrix} U \\ V \end{bmatrix} \sim N_{p+q}\left(\begin{bmatrix} 0 \\ \mu_V \end{bmatrix}, \begin{bmatrix} Q & A^T \\ A & S \end{bmatrix}\right), \begin{bmatrix} Z \\ V \end{bmatrix} \sim N_{r+q}\left(\begin{bmatrix} \mu_Z \\ \mu_V \end{bmatrix}, \begin{bmatrix} R & C^T \\ C & S \end{bmatrix}\right).$$

We allow $Q$ and $R$, but not $S$, to be degenerate. According to a familiar result in probability, the conditional distribution of $U$ given $V = v$ is:

$$U|v \sim N_p(A^T S^{-1}(v - \mu_V), Q - A^T S^{-1} A). \tag{29}$$

We shall study the following situation. Suppose $\lambda \in L(R^p \otimes R^p; R^r)$ and $\theta \in L(R^p \otimes R^p; R^q)$ are symmetric bilinear forms, and $X$ and $Y$ are random vectors of the form:

$$X \equiv Z + \lambda(U \otimes U), Y \equiv V + \theta(U \otimes U),$$

for $Z$, $U$ and $V$ as above. Assume moreover that each of $\|Q\|$, $\|R\|$, $\|S\|$ is $O(\gamma^2)$, for some small number $\gamma$.

### 4.1.a    Definition

We shall say that $W$ approximates $E[X|Y]$ up to $O(\gamma^4)$ if

$$E[h(Y - E[Y]) \otimes (W - X)] = O(\gamma^4), \tag{30}$$

for every $h \in C^1(R^q; R^r)$ with $\max\{\sup_y \|h(y)\|, \sup_y \|Dh(y)\|\} = 1$.

## 4.2    Proposition

*The conditional expectation $E[X|Y]$ is approximated up to $O(\gamma^4)$ (in the sense of (30)) by*

$$W = E[X] + G\hat{Y} + \rho(\hat{Y} \otimes \hat{Y}) - E[\rho(\hat{Y} \otimes \hat{Y})] \in R^p, \tag{31}$$

*where $\hat{Y} \equiv Y - E[Y]$, and the coefficients $G \in L(R^q; R^r)$ and $\rho \in L(R^q \otimes R^q; R^r)$ are given by*

$$G \equiv C^T S^{-1}, \rho(y \otimes y) \equiv (\lambda - G\theta)(A^T S^{-1} y \otimes A^T S^{-1} y). \tag{32}$$

*In the case where* $Z = U$, $\mathrm{Var}(X|Y)$ *is approximated up to* $O(\gamma^4)$ *by*

$$R - C^T S^{-1} C = \mathrm{Var}(Z|V). \qquad (33)$$

**Proof:** Without loss of generality, we may take $\mu_V = 0$ and $\mu_Z = 0$, so $\hat{Y} \equiv Y - \theta(Q)$. The strategy of the proof will be to expand $h(\hat{Y})$ in a Taylor series as
$h(V) + Dh(V)(\theta(U \otimes U - Q)) + \ldots$, so that we can evaluate the left side of (30) as:

$$\mathrm{E}[\mathrm{E}[\{h(V) + Dh(V)(\theta(U \otimes U - Q)) + \ldots\} \otimes (W - X)|V]]$$

$$\mathrm{E}[h(V)\mathrm{E}[(W - X)|V]] + \mathrm{E}[Dh(V)(\theta(U \otimes U - Q))\mathrm{E}[(W - X)|V]] + \ldots. \qquad (34)$$

**Step I.** First we evaluate $\mathrm{E}[W - X|V]$. Applying (29),

$$\mathrm{E}[W - X|V] = \mathrm{E}[X] - \mathrm{E}[\rho(\hat{Y} \otimes \hat{Y})] + \mathrm{E}[G\hat{Y} + \rho(\hat{Y} \otimes \hat{Y}) - Z - \lambda(U \otimes U)|V]$$

$$= (\lambda - G\theta)(Q) - \mathrm{E}[\rho(V \otimes V)] + (G - C^T S^{-1})V$$

$$+ \mathrm{E}[(G\theta - \lambda)(U \otimes U) + \rho([V + \theta(U \otimes U - Q)] \otimes [V + \theta(U \otimes U - Q)])|V]$$

$$= (G - C^T S^{-1})V + \rho(V \otimes V - S) - (\lambda - G\theta)(\mathrm{E}[U \otimes U|V] - Q)$$

$$+ 2\rho(V \otimes \theta\{\mathrm{E}[U \otimes U|V] - Q\}) + f_0^4(V),$$

where $f_i^4(V)$ refers to a linear function of a quadruple tensor product of normal random variables, conditioned on $V$. Using (29) again, we see that

$$\mathrm{E}[U \otimes U|V] = \mathrm{E}[(U - A^T S^{-1} V) \otimes U|V] + A^T S^{-1} V \otimes \mathrm{E}[U|V]$$

$$= Q - A^T S^{-1} A + A^T S^{-1} V \otimes A^T S^{-1} V.$$

Inserting the values of $G$ and $\rho$ specified in (32), and taking $H \equiv A^T S^{-1}$,

$$\mathrm{E}[W - X|V] = 2\rho(V \otimes \theta\{HV \otimes HV - \mathrm{E}[HV \otimes HV]\}) + f_0^4(V). \qquad (35)$$

**Step II.** Calculating in a similar way,

$$\mathrm{E}[\theta(U \otimes U - Q) \otimes (W - X)|V] = \mathrm{Cov}(\theta(U \otimes U), G\hat{Y} + \rho(\hat{Y} \otimes \hat{Y}) - Z - \lambda(U \otimes U)|V)$$

$$= -\mathrm{Cov}(\theta(U \otimes U), Z|V) + f_1^4(V)$$

$$= -\mathrm{Cov}(\theta(U \otimes U), \mathrm{E}[Z|U]|V) - \mathrm{E}[\mathrm{Cov}(\theta(U \otimes U), Z|U)|V] + f_1^4(V).$$

The last line follows from the identity "$\mathrm{Var}(X) = \mathrm{Var}(\mathrm{E}[X|U]) + \mathrm{E}[\mathrm{Var}(X|U)]$". The conditional covariance in the second term is zero. By (29) $\mathrm{E}[Z|U]$ is a linear function of $U$, and since third moments of $U$ vanish by the normality assumption, we conclude that

$$E[\theta(U \otimes U - Q) \otimes (W - X) | V] = f_1^4(V).$$

**Step III.** Let $h$ be as in Definition 4.1.a. Inserting the results of Steps I and II into (34),

$$E[h(\hat{Y}) \otimes (W - X)] =$$

$$E\left[h(V) \otimes \{2\rho(V \otimes \theta\{HV \otimes HV - E[HV \otimes HV]\}) + f_0^4(V)\} + Dh(V) \otimes f_1^4(V) + f_2^4(V)\right]$$

and this is $O(\gamma^4)$ because $h(V) = h(0) + Dh(0)V + O(\gamma^2)$, and third moments of $V$ vanish by the normality assumption. Thus we have verified that $W$ approximates $E[X|Y]$ up to $O(\gamma^4)$.

**Step IV.** It remains to check the covariance assertion (33). By analogy with Step II,

$$E[\theta(U \otimes U - Q) \otimes (W - X) \otimes (W - X) | V] = E[\theta(U \otimes U - Q) \otimes Z \otimes Z | V] + f_3^4(V) =$$

$$E[\theta(U \otimes U - Q) \otimes \text{Var}(Z|U) | V] + E[\theta(U \otimes U - E[U \otimes U | V]) \otimes E[Z|U] \otimes E[Z|U] | V]$$

$$+ f_3^4(V)$$

which follows from a generalization of "$\text{Var}(X) = \text{Var}(E[X|U]) + E[\text{Var}(X|U)]$". However $E[Z|U] = U$, and $\text{Var}(Z|U) = 0$ under the assumption that $Z = U$. Thus

$$E[\theta(U \otimes U - Q) \otimes (W - X) \otimes (W - X) | V] = f_4^4(V). \quad (36)$$

We also have

$$E[(W - X) \otimes (W - X) | V] = \text{Var}(G\hat{Y} + \rho(\hat{Y} \otimes \hat{Y}) - Z - \lambda(U \otimes U) | V)$$

$$= \text{Var}(-Z|V) + f_5^4(V)$$

$$= R - C^T S^{-1} C + f_5^4(V). \quad (37)$$

In view of the results of Step III, it suffices to show that

$$E\left[h(\hat{Y}) \otimes \{(W - X) \otimes (W - X) - (R - C^T S^{-1} C)\}\right] = O(\gamma^4).$$

We substitute the Taylor expansion for $h(\hat{Y})$, condition on $V$ inside the expectation, and apply (36) and (37), to obtain (33). ◊

## 5    Intrinsic Approximation of a Diffusion Process

### 5.1    Geometry Induced by a Diffusion Process

Consider a (possibly degenerate) Markov diffusion process $\{X_t, 0 \le t \le \delta\}$ on $N \cong R^p$, written in local coördinates as

$$dX_t^i = b^i(X_t)\,dt + \sum_{j=1}^{p} \sigma_j^i(X_t)\,dW_t^j, \ i = 1, 2, ..., p, \tag{38}$$

where $\sum b^i \frac{\partial}{\partial x_i}$ is a vector field on $R^p$, $\sigma(x) \equiv (\sigma_j^i(x)) \in L(R^p; T_x R^p)$, and $W$ is a Wiener process in $R^p$. We assume for simplicity that the coefficients $b^i$, $\sigma_j^i$ are $C^3$ with bounded first derivative. Such a $\sigma$ induces a $C^2$ semi-definite metric $\langle .|.\rangle$ on the cotangent bundle, which we call the **diffusion variance semi-definite metric**, by the formula

$$\langle dx^i | dx^k \rangle_x \equiv (\sigma \cdot \sigma(x))^{ik} \equiv \sum_{j=1}^{p} \sigma_j^i(x)\,\sigma_j^k(x) . \tag{39}$$

This semi-definite metric is actually intrinsic: changing coördinates for the diffusion will give a different matrix $(\sigma_j^i)$, but the same semi-definite metric. The $p \times p$ matrix $((\sigma \cdot \sigma)^{ij})$ defined above induces a linear transformation $\alpha(x): T_x^* N \to T_x N$, i.e. from the cotangent space to the tangent space at $x$, namely

$$\alpha(x)(dx^i) \equiv \sum (\sigma \cdot \sigma)^{ij} \partial/\partial x_j . \tag{40}$$

Let us make a **constant-rank assumption**, i.e. that there exists a rank $r$ vector bundle $E \to N$, a sub-bundle of the tangent bundle, such that $E_x = \text{range}(\sigma(x)) \subseteq T_x N$ for all $x \in N$. Darling [2] presents a global geometric construction of a **canonical sub-Riemannian connection** $\nabla^\circ$ for $\langle .|.\rangle$, with respect to a generalized inverse $g$, i.e. a vector bundle isomorphism $g: TN \to T^*N$ such that

$$\alpha(x) \bullet g(x) \bullet \alpha(x) = \alpha(x) .$$

In local coördinates, $g(x)$ is expressed by a Riemannian metric tensor $(g_{rs})$, such that if $\alpha^{ij} \equiv (\sigma \cdot \sigma)^{ij}$, then

$$\sum_{r,s} \alpha^{ir} g_{rs} \alpha^{sj} = \alpha^{ij} . \tag{41}$$

The local connector $\Gamma(x) \in L(T_x R^p \otimes T_x R^p; T_x R^p)$ for $\nabla^\circ$ can be written in the form:

$$2g(\Gamma(x)(u \otimes v)) \cdot w = D\langle g(v)|g(w)\rangle(u) + D\langle g(w)|g(u)\rangle(v) - D\langle g(u)|g(v)\rangle(w), \tag{42}$$

where $g(\Gamma(x)(u \otimes v))$ is a 1-form, acting on the tangent vector $w$. This formula coincides with the formula for the Levi-Civita connection in the case where $\langle .|.\rangle$ is non-degenerate; for more details, see Darling [2].

### 5.2   Intrinsic Description of the Process

The **intrinsic** version of (38) is to describe $X$ as a diffusion process on the manifold $N$ with generator

$$L \equiv \xi + \frac{1}{2}\Delta \tag{43}$$

where $\Delta$ is the (possibly degenerate) Laplace-Beltrami operator associated with the diffusion variance, and $\xi$ is a vector field, whose expressions in the local coördinate system $\{x^1, ..., x^p\}$ are as follows:

$$\Delta = \sum_{i,j} (\sigma \cdot \sigma)^{ij} \{D_{ij} - \sum_k \Gamma_{ij}^k D_k\}, \quad \xi = \sum_k \{b^k + \frac{1}{2} \sum_{i,j} (\sigma \cdot \sigma)^{ij} \Gamma_{ij}^k\} D_k, \qquad (44)$$

where the Christoffel symbols $\{\Gamma_{ij}^k\}$ are related to the conector by

$$\Gamma_{ij}^k \equiv dx^k \cdot \Gamma(x) \left(\frac{\partial}{\partial x_i} \otimes \frac{\partial}{\partial x_j}\right).$$

Uncertainty about the initial value will be expressed by writing

$$X_0 \equiv \exp_{x_0}(U_0), \qquad (45)$$

where $U_0$ is a random variable in $T_{x_0}N$, with covariance tensor $\Sigma_0 \in T_xN \otimes T_xN$, and the exponential map is taken with respect to the connection defined in (42).

## 5.3   Intrinsic Linearization About the Deterministic Solution

Let $\{x_t, 0 \leq t \leq \delta\}$ be the solution of the ODE associated with the vector field $\xi$, started at the same $x_0$ as is mentioned in (45). In other words, $x_t = \phi_t(x)$, for $0 \leq t \leq \delta$, where $\{\phi_t, 0 \leq t \leq \delta\}$ is the flow of the vector field $\xi$ on $N$. We may compute the following intrinsic quantities:

$$\tau_s^t \equiv D(\phi_t \bullet \phi_s^{-1})(x_s) \in L(T_{x_s}N; T_{x_t}N); \qquad (46)$$

$$\Pi_t \equiv \Sigma_0 + \int_0^t (\phi_{-s})_* \langle .|.\rangle_{x_s} ds = \Sigma_0 + \int_0^t \tau_s^0 (\sigma \cdot \sigma)(x_s)(\tau_s^0)^T ds \in T_{x_0}N \otimes T_{x_0}N; \qquad (47)$$

$$\Xi_t \equiv (\phi_t)_* \Pi_t = \int_0^t \tau_s^t (\sigma \cdot \sigma)(x_s)(\tau_s^t)^T ds + \tau_0^t \Sigma_0 (\tau_0^t)^T \in T_{x_t}N \otimes T_{x_t}N. \qquad (48)$$

## 5.4   Approximate Intrinsic Location Parameter

We recall from Darling [2] that there exists a vector in the tangent space $T_{\psi(x_\delta)}M$ which supplies a coördinate-independent replacement for the notion of expected value $\psi(X_\delta)$. This vector, denoted $I_{x_0, \Sigma_0}[\psi(X_\delta)]$, is called the **approximate intrinsic location parameter (AILP)** of $\psi(X_\delta)$ in the tangent space $T_{\psi(x_\delta)}M$. We here omit any discussion of how the AILP is derived from the study of manifold-valued martingales, or its relation to harmonic mappings, but merely state the formula

$$I_{x_0, \Sigma_0}[\psi(X_\delta)] = \frac{1}{2}\{\nabla d\psi(x_\delta)(\Xi_\delta) + \psi_*\left[\nabla d\phi_\delta(x_0)(\Pi_\delta) - \int_0^\delta \tau_t^\delta \nabla d\phi_t(x_0)(d\Pi_t)\right]\}. \qquad (49)$$

Here $\nabla d\psi$ denotes the second fundamental form of the map $\psi:N \to M$; if $N$ has connector $\Gamma(.)$ and $M$ has connector $\bar{\Gamma}(.)$, the formula for $\nabla d\psi(x)(v \otimes w)$ can be written as

$$D^2\psi(x)\,(v\otimes w) - D\psi(x)\,\Gamma(x)\,(v\otimes w) + \bar{\Gamma}(y)\,(D\psi(x)\,v\otimes D\psi(x)\,w)\,. \tag{50}$$

In the particular case where $\psi$ is the identity, we obtain $I_{x_0,\Sigma_0}[X_\delta]$, namely the AILP of $X_\delta$ in the tangent space $T_{x_0}N$, namely

$$I_{x_0,\Sigma_0}[X_\delta] = \frac{1}{2}\{\nabla d\phi_\delta(x_0)\,(\Pi_\delta) - \int_0^\delta \tau_t^\delta \nabla d\phi_t(x_0)\,(d\Pi_t)\}\,. \tag{51}$$

## 5.5 A One-parameter Family of Diffusion Processes

Consider a family of $N$-valued diffusion processes $\{X^\varepsilon, \varepsilon \geq 0\}$ on the time interval $[0,\delta]$, where $X^\varepsilon$ has generator $\xi + \varepsilon^2 \Delta/2$, realized in local coördinates by

$$X_t^\varepsilon = \exp_{x_0}(\varepsilon U_0) + \int_0^t \xi(X_s^\varepsilon)\,ds - \varepsilon^2 \int_0^t \zeta(X_s^\varepsilon)\,ds + \int_0^t \varepsilon\sigma(X_s^\varepsilon)\,dW_s\,. \tag{52}$$

We use the notation

$$\zeta(x) \equiv \frac{1}{2}\Gamma(x)\,(\sigma(x)\cdot\sigma(x))\,.$$

In the case $\varepsilon = 0$, the solution is deterministic, namely $\{x_t, 0 \leq t \leq \delta\}$. Differentiating with respect to $\varepsilon$, and using Lemma 2.5, gives formulas for

$$\Lambda_t \equiv \left.\frac{\partial X_t^\varepsilon}{\partial \varepsilon}\right|_{\varepsilon=0},\ \Lambda_{2,t} \equiv \left.\frac{\partial^2 X_t^\varepsilon}{\partial \varepsilon^2}\right|_{\varepsilon=0},$$

namely (see Darling [2])

$$\Lambda_t = \tau_0^t U_0 + \tau_0^t \int_0^t \tau_s^0 \sigma(x_s)\,dW_s = \tau_0^t U_0 + \int_0^t \tau_s^t \sigma(x_s)\,dW_s, \tag{53}$$

$$\Lambda_{2,t} = -\tau_0^t \Gamma(x_0)\,(U_0 \otimes U_0) + \int_0^t \tau_s^t \{[D^2\xi(x_s)\,(\Lambda_s \otimes \Lambda_s) - 2\zeta(x_s)]\,ds + 2D\sigma(x_s)\,(\Lambda_s)\,dW_s\}\,, \tag{54}$$

where $\tau_s^t$ is as given in (46). A natural object is the process $\{\Omega_t, 0 \leq t \leq \delta\}$ given by

$$\Omega_t \equiv (\tau_t^0 \Lambda_t) \otimes (\tau_t^0 \Lambda_t) = \left(U_0 + \int_0^t \tau_s^0 \sigma(x_s)\,dW_s\right) \otimes \left(U_0 + \int_0^t \tau_s^0 \sigma(x_s)\,dW_s\right) \in T_{x_0}N \otimes T_{x_0}N\,. \tag{55}$$

With reference to (47) and (48), taking $\Psi_t \equiv \Lambda_t \otimes \Lambda_t$, we have

$$\mathrm{Var}(\tau_t^0 \Lambda_t) = E[\Omega_t] = \Pi_t,\ \mathrm{Var}(\Lambda_t) = E[\Psi_t] = \Xi_t. \tag{56}$$

## 5.6 Proposition

*Let $u(\varepsilon) \equiv \exp_{x_\delta}^{-1}(X_\delta^\varepsilon)$. Then $u'(0) = \Lambda_\delta$, and*

$$u''(0) = \nabla d\phi_\delta(x_0)(\Omega_\delta) - \int_0^\delta \tau_t^\delta \nabla d\phi_t(x_0)(d\Pi_t) + H_\delta, \tag{57}$$

in the notation of (46), (47), and (55), where

$$H_\delta \equiv 2\int_0^\delta \tau_t^\delta \{\nabla_{\Lambda_t}\sigma(x_t)\,dW_t - \nabla d\phi_t(x_0)(\tau_t^0 \Lambda_t \otimes \tau_t^0 \sigma(x_t)\,dW_t)\}. \tag{58}$$

Moreover

$$E\left[u'(0) + \frac{1}{2}u''(0)\right] = I_{x_0,\Sigma_0}[X_\delta]. \tag{59}$$

**Proof:** It follows immediately from Corollary 2.6 that $u'(0) = \Lambda_\delta$ and

$$u''(0) = \Lambda_{2,\delta} + \Gamma(x_\delta)(\Lambda_\delta \otimes \Lambda_\delta).$$

Now $\Psi_t \equiv \Lambda_t \otimes \Lambda_t$ satisfies

$$d\Psi_t = \{D\xi(x_t)\Psi_t + \Psi_t(D\xi(x_t))^T + (\sigma \cdot \sigma)(x_t)\}\,dt + \sigma(x_t)\,dW_t \otimes \Lambda_t + \Lambda_t \otimes \sigma(x_t)\,dW_t.$$

Copying the structure of calculations presented at the end of the proof of the main result of Darling [2], it follows that

$$d\{(\phi_{\delta-t})_*\nabla d\phi_t(x_0)(\Omega_t)\} - (\phi_{\delta-t})_*\nabla d\phi_t(x_0)\,d\Omega_t = \frac{\partial}{\partial t}\{(\phi_{\delta-t})_*\nabla d\phi_t(x_0)\}(\Omega_t)\,dt$$

$$= \left\{\left[\tau_t^\delta D^2\xi(x_t) + \frac{\partial}{\partial t}\{\tau_t^\delta \Gamma(x_t)\}\right](\Psi_t) + \tau_t^\delta \Gamma(x_t)[D\xi(x_t)\Psi_t + \Psi_t(D\xi(x_t))^T]\right\}dt$$

$$= \left[\tau_t^\delta D^2\xi(x_t) + \frac{\partial}{\partial t}\{\tau_t^\delta \Gamma(x_t)\}\right](\Psi_t)$$

$$+ \tau_t^\delta \Gamma(x_t)[d\Psi_t - \sigma(x_t)\,dW_t \otimes \Lambda_t - \Lambda_t \otimes \sigma(x_t)\,dW_t - (\sigma \cdot \sigma)(x_t)\,dt]$$

$$= \tau_t^\delta[D^2\xi(x_t)(\Psi_t) - \tau_t^\delta \Gamma(x_t)(\sigma \cdot \sigma(x_t))]\,dt + d(\tau_t^\delta \Gamma(x_t)(\Psi_t)) - 2\tau_t^\delta \Gamma(x_t)(\Lambda_t \otimes \sigma(x_t)\,dW_t).$$

Since $\nabla d\phi_0 = 0$, it follows upon integration from 0 to $\delta$ that in $T_{x_\delta}M$,

$$\nabla d\phi_\delta(x_0)(\Omega_\delta) - \int_0^\delta (\phi_{\delta-t})_*(\nabla d\phi_t(x_0))\,d\Omega_t =$$

$$\int_0^\delta \tau_t^\delta[D^2\xi(x_t)(\Psi_t) - \Gamma(x_t)(\sigma \cdot \sigma(x_t))]\,dt + \Gamma(x_\delta)(\Psi_\delta) - \tau_0^\delta \Gamma(x_0)(\Psi_0)$$

$$- 2\int_0^\delta \tau_t^\delta \Gamma(x_t)(\Lambda_t \otimes \sigma(x_t)\,dW_t).$$

Since $(\phi_{\delta-t})_* = \tau_t^\delta$, we have

$$\Lambda_{2,\delta} + \Gamma(x_\delta)(\Lambda_\delta \otimes \Lambda_\delta) = \quad (60)$$

$$\int_0^\delta \tau_s^\delta \{[D^2\xi(x_s)(\Psi_s) - 2\zeta(x_s)]ds + 2D\sigma(x_s)(\Lambda_s)dW_s\} + \Gamma(x_\delta)(\Psi_\delta) - \tau_0^\delta \Gamma(x_0)(\Psi_0)$$

$$= 2\int_0^\delta \tau_t^\delta \{\Gamma(x_t)(\Lambda_t \otimes \sigma(x_t)dW_t) + D\sigma(x_t)(\Lambda_t)dW_t\} + \nabla d\phi_\delta(x_0)(\Omega_\delta) - \int_0^\delta \tau_t^\delta \nabla d\phi_t(x_0)d\Omega_t.$$

However from (55),

$$d\Omega_t = 2(\overset{0}{\tau}_t \Lambda_t) \otimes \overset{0}{\tau}_t \sigma(x_t)dW_t + \overset{0}{\tau}_t(\sigma \cdot \sigma)(x_t)(\overset{0}{\tau}_t)^T dt,$$

and for $e \in R^p$, (1) implies that

$$\Gamma(x_t)(\Lambda_t \otimes \sigma(x_t)e) + D\sigma(x_t)(\Lambda_t)e = \nabla_{\Lambda_t}\sigma(x_t)e.$$

and (57) follows. Note that $E[\nabla d\phi_\delta(x_0)(\Omega_\delta)] = \nabla d\phi_\delta(x_0)(\Pi_\delta)$, and so (59) follows from (51) and (57).                                                                                  ◊

## 6 Intrinsic Estimation of a Diffusion Using a Single Observation

### 6.1 Observation of the Process

We are given a $C^3$ function $\psi: N \to M$, where $M$ is a Riemannian manifold of dimension $q$. Let $\beta(y) \in T_y M \otimes T_y M$ be the inverse metric tensor at $y \in M$, which can be interpreted as the covariance of a random vector in $T_y M$. Consider a single observation $Y_1$ of the form:

$$Y_1 \equiv \exp_{\psi(X_\delta)} V_1 \in M,$$

where $V_1$ is a mean-zero random vector in $T_{\psi(X_\delta)}M$, whose covariance is $\beta(y)$ when $\psi(X_\delta) = y$, but which is otherwise independent of $U_0$ and the Wiener process $W$.

### 6.2 Orders of Magnitude of Noise Terms

We shall suppose that, for some small number $\gamma$, the tensor fields $\alpha \equiv \sigma \cdot \sigma$ (see (40)) and $\beta$ (see Section 6.1) satisfy

$$\alpha(x_t) = \gamma^2 \alpha_0(x_t), \ 0 \le t \le \delta; \ \beta(\psi(x_\delta)) = \gamma^2 \beta_0(\psi(x_\delta)); \quad (61)$$

where $\alpha_0$ is some other semi-definite metric, and $\beta_0$ another metric. Also assume that, with respect to the metric $g$ appearing in (41), the distribution of $U_0 \equiv \exp_{x_0}^{-1}(X_0)$ satisfies:

$$\|E[U_0]\| = O(\gamma^4), \ \|\Sigma_0\| \equiv \|\text{Var}(U_0)\| = O(\gamma^2), \ \|E[T(U_0, U_0, U_0)]\| = O(\gamma^4), \quad (62)$$

for arbitrary tensor fields $T$ of type $(1,3)$ of norm 1.

### 6.3 Theorem

Consider the random vector $U_\delta \oplus Z_\delta \in T_{x_\delta} N \oplus T_{\psi(x_\delta)} M$ given by

$$U_\delta \equiv \exp^{-1}_{x_\delta}(X_\delta) , \quad Z_\delta \equiv \exp^{-1}_{\psi(x_\delta)}(Y_1) . \tag{63}$$

(i) Under the assumptions (61) and (62), the joint distribution of $U_\delta$ and $Z_\delta$ satisfies

$$E\begin{bmatrix} U_\delta \\ Z_\delta \end{bmatrix} = \begin{bmatrix} I_{x_0, \Sigma_0}[X_\delta] \\ I_{x_0, \Sigma_0}[\psi(X_\delta)] \end{bmatrix} + O(\gamma^4) ; \tag{64}$$

where the expressions in (64) are approximate intrinsic location parameters (see Section 5.4), $\Xi_\delta$ is given by (48), and $J \equiv D\psi(x_\delta)$.

(ii) $E[U_\delta | Z_\delta]$ is approximated up to $O(\gamma^4)$ (in the sense of (30)) by

$$I_{x_0, \Sigma_0}[X_\delta] + G\hat{Z}_\delta + \rho(\hat{Z}_\delta \otimes \hat{Z}_\delta) - E[\rho(\hat{Z}_\delta \otimes \hat{Z}_\delta)] , \tag{65}$$

where $\hat{Z}_\delta \equiv Z_\delta - I_{x_0, \Sigma_0}[\psi(X_\delta)]$, and $G \in L(T_{y_\delta} M; T_{x_\delta} N)$ is analogous to the Kalman gain, namely

$$G \equiv \Xi_\delta J^T [J\Xi_\delta J^T + \beta(\psi(x_\delta))]^{-1} , \tag{66}$$

$$\rho(z \otimes z) \equiv \frac{1}{2} \{ [I - GJ] \nabla d\phi_\delta(x_0) (\tau^0_\delta Gz \otimes \tau^0_\delta Gz) - G\nabla d\psi(x_\delta)(Gz \otimes Gz) \} , \tag{67}$$

$$E[\rho(\hat{Z}_\delta \otimes \hat{Z}_\delta)] = \rho(GJ\Xi_\delta) + O(\gamma^4) .$$

(iii) $\mathrm{Var}(U_\delta | Z_\delta)$ is approximated up to $O(\gamma^4)$ by $(I - GJ)\Xi_\delta$.

(iv) If $\hat{U}_\delta$ denotes the difference between $U_\delta$ and (65), and if $T$ is a tensor field of type $(1,3)$ on $N$ of norm 1, then

$$\left\| E[T(\hat{U}_\delta, \hat{U}_\delta, \hat{U}_\delta)] \right\| = O(\gamma^4) . \tag{68}$$

**Proof: Step I.** We continue to study the family of $N$-valued diffusion processes $\{X^\varepsilon, \varepsilon \geq 0\}$ given by (52), whose first and second derivatives with respect to $\varepsilon$ were presented in (53) and (54). As for the third derivative,

$$\Lambda_{3,t} \equiv \left. \frac{\partial^3 X^\varepsilon_t}{\partial \varepsilon^3} \right|_{\varepsilon = 0} = \Theta^U_0 + \Theta^W_t + \int_0^t D\xi(x_s)(\Lambda_{3,s}) ds$$

where $\Theta^U_0$ is a trilinear function of $U_0$, and $\Theta^W_t$ is a sum of integrals whose expectation is zero at $\varepsilon = 0$, and whose variance is $O(\gamma^4)$ or higher order. We shall use the same symbols from line to line, even though we may add or delete terms of each type as we go along. We obtain

$$\Lambda_{3,t} = \tau_0^t \Theta_0^U + \int_0^t \tau_s^t d\Theta_s^W,$$

which can be condensed to

$$\Lambda_{3,\delta} = \Theta_0^U + \Theta_\delta^W. \tag{69}$$

**Step II.** Define a random path $\{x(\varepsilon), 0 \leq \varepsilon \leq 1\}$ in $N$ by $x(\varepsilon) \equiv X_\delta^\varepsilon$, and a random path $\{y(\varepsilon), 0 \leq \varepsilon \leq 1\}$ in $M$ by $y(\varepsilon) \equiv \psi(x(\varepsilon))$. Both are thrice differentiable in $\varepsilon$, and the derivatives of $y$ at $\varepsilon = 0$ are given as follows.

$$y'(0) = J\Lambda_\delta \sim N_q(0, J\Xi_\delta J^T), \tag{70}$$

for $J \equiv D\psi(x_\delta)$ and $\Xi_\delta$ as in (48). With a further differentiation, and use of (1) and (50), we have

$$\frac{\nabla y'}{\partial \varepsilon}(0) = \nabla d\psi(x_\delta)(\Lambda_\delta \otimes \Lambda_\delta) + Ju''(0) \in T_{\psi(x_\delta)}M,$$

for $u(\varepsilon) \equiv \exp_{x_\delta}^{-1}(X_\delta^\varepsilon)$ as in Proposition 5.6. Performing a further covariant differentiation, we find from (53), (54), and (69) that

$$\frac{\nabla^2 y'}{\partial \varepsilon^2}(0) = \Theta_0^U + \Theta_\delta^W + O(\gamma^4).$$

In particular,

$$E\left[\frac{\nabla^2 y'}{\partial \varepsilon^2}(0)\right] = O(\gamma^4).$$

**Step III.** Let $F$ be a smooth section of $\text{Hom}(R^q; TM)$ (i.e. $F(y) \in L(R^q; T_y M)$ for all $y \in M$) such that $F(y)F(y)^T = \beta(y)$, the observation covariance metric, and define a random vector field along $y$ by

$$V(\varepsilon) \equiv \varepsilon F(y(\varepsilon))Z \in T_{y(\varepsilon)}M, \tag{71}$$

where $Z$ is a $N_q(0, I)$ random variable in $R^q$, independent of $U_0$ and $W$. Our object of interest is the random variable

$$\zeta(\varepsilon) \equiv \exp_{y(0)}^{-1}(\exp_{y(\varepsilon)} V(\varepsilon)) \in T_{y(0)}M,$$

since $\zeta(1) = Z_\delta$. By (71), $V(0) = 0$, and

$$\frac{\nabla V}{\partial \varepsilon} = F(y)Z + \varepsilon \frac{\nabla F}{\partial \varepsilon}Z, \quad \frac{\nabla^2 V}{\partial \varepsilon^2} = 2\frac{\nabla F}{\partial \varepsilon}Z + \varepsilon \frac{\nabla^2 F}{\partial \varepsilon^2}Z, \quad \frac{\nabla^3 V}{\partial \varepsilon^3} = 3\frac{\nabla^2 F}{\partial \varepsilon^2}Z + \varepsilon \frac{\nabla^3 F}{\partial \varepsilon^3}Z, \tag{72}$$

noting that, by (1), $\frac{\nabla F}{\partial \varepsilon}Z = DF(y)(y')Z + \Gamma(y)(F(y)Z \otimes y')$.

According to Proposition 2.4, and the formulas of Step II, writing $y_\delta$ for $y(0) \equiv \psi(x_\delta)$,

$$\zeta'(0) = J\Lambda_\delta + F(y_\delta) Z + O(\gamma^4) ;$$

$$\zeta''(0) = \nabla d\psi(x_\delta) (\Lambda_\delta \otimes \Lambda_\delta) + Ju''(0) + 2\frac{\nabla F}{\partial \varepsilon}(y_\delta) Z + O(\gamma^4) ;$$

$$\zeta^{(3)}(0) = \Theta_0^U + \Theta_\delta^W + 3\frac{\nabla^2 F}{\partial \varepsilon^2}(y_\delta) Z + R(y_\delta) (J\Lambda_\delta, F(y_\delta) Z) (J\Lambda_\delta + 2F(y_\delta) Z) + O(\gamma^4) .$$

Taking a third order Taylor approximation at $\varepsilon = 0$, with $\zeta(0) = 0$,

$$Z_\delta \equiv \zeta(1) = J\{\Lambda_\delta + \frac{1}{2}u''(0)\} + \frac{1}{2}\nabla d\psi(x_\delta) (\Lambda_\delta \otimes \Lambda_\delta) + \Theta_0^U + \Theta_\delta^W \tag{73}$$

$$+ \left\{F + \frac{\nabla F}{\partial \varepsilon} + \frac{1}{2}\frac{\nabla^2 F}{\partial \varepsilon^2}\right\}(y_\delta) Z + \frac{1}{6}R(y_\delta) (J\Lambda_\delta, F(y_\delta) Z) (J\Lambda_\delta + 2F(y_\delta) Z) + O(\gamma^4) .$$

In the special case where $\psi$ is the identity, and $Z$ is suppressed, we have the analogous formula

$$U_\delta = \Lambda_\delta + \frac{1}{2}u''(0) + \Theta_0^U + \Theta_\delta^W + O(\gamma^4) . \tag{74}$$

**Step IV.** We see from Proposition 5.6 that (74) can be written in the form

$$U_\delta = I_{x_0, \Sigma_0}[X_\delta] - E[\lambda(\Lambda_\delta \otimes \Lambda_\delta)] + \Lambda_\delta + \lambda(\Lambda_\delta \otimes \Lambda_\delta) + h_\delta^U, \tag{75}$$

where the bilinear mapping $\lambda$ and the random vector $h_\delta^U$ are given, respectively, by

$$\lambda(w \otimes w) \equiv \frac{1}{2}\nabla d\phi_\delta(x_0) (\tau_\delta^0 w \otimes \tau_\delta^0 w) ,$$

$$h_\delta^U \equiv 2\int_0^\delta \tau_t^\delta \{\nabla_{\Lambda_t}\sigma(x_t) dW_t - \nabla d\phi_t(x_0) (\tau_t^0 \Lambda_t \otimes \tau_t^0 \sigma(x_t) dW_t)\} + \Theta_0^U + \Theta_\delta^W + O(\gamma^4) .$$

Our assumptions in Section 6.2 imply that $E[h_\delta^U] = O(\gamma^4)$ and $\text{Var}(h_\delta^U) = O(\gamma^4)$. It also follows from (49), (59), (56), and (73) that

$$Z_\delta = I_{x_0, \Sigma_0}[\psi(X_\delta)] - E[\theta(\Lambda_\delta \otimes \Lambda_\delta)] + J\Lambda_\delta + F(y_\delta) Z + \theta(\Lambda_\delta \otimes \Lambda_\delta) + h_\delta^Z, \tag{76}$$

where $\theta$ is the bilinear mapping given by

$$\theta(w \otimes w) \equiv \frac{1}{2}\nabla d\psi(x_\delta) (w \otimes w) + J\lambda(w \otimes w) ,$$

and likewise $E[h_\delta^Z] = O(\gamma^4)$ and $\text{Var}(h_\delta^Z) = O(\gamma^4)$. This depends in particular on the fact that $\left\|\frac{\nabla V}{\partial \varepsilon}\right\|$ is $O(\gamma)$ at $\varepsilon = 0$, and so

$$\text{Var}\left(\frac{\nabla F}{\partial \varepsilon}Z\right) = O(\gamma^4) .$$

The first moment formula (64) follows immediately from (75) and (76). We are now in a position to apply Proposition 4.2. Observe that

$$\begin{bmatrix} \Lambda_\delta \\ J\Lambda_\delta + F(y_\delta)Z \end{bmatrix} \sim N_{p+q}\left( \begin{bmatrix} 0 \\ 0 \end{bmatrix}, \begin{bmatrix} \Xi_\delta & \Xi_\delta J^T \\ J\Xi_\delta & J\Xi_\delta J^T + \beta(y_\delta) \end{bmatrix} \right), \quad (77)$$

and the matrices $C$ and $A$ of Proposition 4.2 are both the same here. According to (31), the conditional expectation $E[U_\delta | Z_\delta]$ is approximated up to $O(\gamma^4)$ (in the sense of (30)) by an expression of the form (65), with $G$ given by (66), and

$$\rho(z \otimes z) = (\lambda - G\theta)(Gz \otimes Gz)$$

$$= \{\lambda - G\left(\frac{1}{2}\nabla d\psi(x_\delta) + J\lambda\right)\}(Gz \otimes Gz),$$

which can be expanded into the form (67). It follows from (77) that, up to $O(\gamma^4)$, $E[\rho(\hat{Z}_\delta \otimes \hat{Z}_\delta)]$ is given by $\rho(GJ\Xi_\delta)$. According to (33) and (77), $Var(U_\delta | Z_\delta)$ is approximated up to $O(\gamma^4)$ by

$$\Xi_\delta - \Xi_\delta J^T (J\Xi_\delta J^T + \beta(y_\delta))^{-1} J\Xi_\delta,$$

which simplifies to $(I - GJ)\Xi_\delta$.

**Step V.** It remains to verify (68). We see from (65) and (75) that $\hat{U}_\delta$ is a random variable whose mean is $O(\gamma^4)$, and which takes the form of a mean-zero Gaussian random variable $S \equiv \Lambda_\delta - G\hat{Z}_\delta$, plus an $O(\gamma^2)$ term. For any tensor field $T$ of type $(1, 3)$ on $N$ of norm 1,

$$E[T(\hat{U}_\delta, \hat{U}_\delta, \hat{U}_\delta)] = E[T(S, S, S)] + O(\gamma^4).$$

However $E[T(S, S, S)] = 0$ since $J$ is mean-zero Gaussian, and (68) follows. ◊

**Acknowledgments:** The author wishes to thank Ofer Zeitouni for perceptive comments on this work, and the Statistics Department at the University of California at Berkeley for its hospitality during the writing of this article.